\documentclass[11pt]{article}
\usepackage[cm]{fullpage}
\usepackage[small]{titlesec}
\usepackage{amscd}
\usepackage{amsmath}
\usepackage{amsfonts}

\newcommand{\ga}{\gamma}
\newcommand{\de}{\partial} 
\newcommand{\db}{\overline{\partial}}

\newcommand{\dd}[1]{\partial_{#1}}

\newcommand{\Ric}{\mathrm{Ric}}

\newcommand{\ov}[1]{\overline{#1}}

\newcommand{\tr}[2]{\textrm{tr}_{#1}{#2}}

\newcommand{\vp}{\varphi}

\newcommand{\ve}{\varepsilon}

\begin{document}
\newcounter{remark}
\newcounter{theor}
\setcounter{remark}{0}
\setcounter{theor}{1}
\newtheorem{claim}{Claim}
\newtheorem{theorem}{Theorem}[section]
\newtheorem{proposition}{Proposition}[section]
\newtheorem{question}{question}[section]
\newtheorem{lemma}{Lemma}[section]
\newtheorem{defn}{Definition}[theor]

\newtheorem{corollary}{Corollary}[section]
\newenvironment{proof}[1][Proof]{\begin{trivlist}
\item[\hskip \labelsep {\bfseries #1}]}{\end{trivlist}}

\newenvironment{remark}[1][Remark]{\addtocounter{remark}{1} \begin{trivlist}
\item[\hskip
\labelsep {\bfseries #1  \thesection.\theremark}]}{\end{trivlist}}

\centerline{\bf \large Estimates for the complex Monge-Amp\`ere equation}
\smallskip

\centerline{\bf \large on Hermitian and balanced manifolds\footnote{Research supported in part by National Science Foundation grant DMS-08-48193.  The second-named author is also supported in part by a Sloan Foundation fellowship.}}

\bigskip
\bigskip
\centerline{\bf Valentino Tosatti and Ben Weinkove}

\bigskip

\begin{abstract} We generalize Yau's estimates for the complex Monge-Amp\`ere equation on compact manifolds in the case when the background metric is no longer K\"ahler.  We prove $C^{\infty}$ \emph{a priori} estimates for a solution of the complex Monge-Amp\`ere equation when the background metric is
Hermitian (in complex dimension two) or balanced (in higher dimensions), giving an alternative proof of a theorem of Cherrier.  We relate this  to recent results of Guan-Li.

\end{abstract}

\bigskip

\section{Introduction}

Let $(M,g)$ be a compact Hermitian manifold of complex dimension $n\ge 2$.  Write $\omega$ for the real $(1,1)$-form
$$\omega = \sqrt{-1} \sum_{i,j} g_{i \ov{j}}  dz^i \wedge d\ov{z^{j}}.$$
We consider the complex Monge-Amp\`ere equation
\begin{equation} \label{ma0}
(\omega + \sqrt{-1} \partial \ov{\partial} \varphi)^n = e^F \omega^n, \quad \textrm{with} \quad \omega+ \sqrt{-1} \partial \ov{\partial} \varphi >0,
\end{equation}
for a real-valued function $\varphi$, where $F$ is a given smooth real-valued function.  Since (\ref{ma0}) is invariant under the addition of constants to $\varphi$, we may assume for convenience that $\varphi$ satisfies the normalization condition,
\begin{equation} \label{sup}
\sup_M \varphi =0.
\end{equation}
When $\omega$ is K\"ahler (that is, $d\omega=0$), Calabi \cite{C} showed that  solutions $\varphi$ to (\ref{ma0}),  (\ref{sup}), are unique.  In fact, the same proof (see Remark \ref{unique}.1 below) shows that this result can be extended for general Hermitian $\omega$. 

In the seminal paper \cite{Y}, Yau showed in the K\"ahler case that one can find a  smooth solution $\varphi$ of (\ref{ma0})  whenever $F$ satisfies the necessary condition
\begin{equation} \label{F}
\int_M e^F \omega^n = \int_M \omega^n.
\end{equation}
Yau proved this result using the continuity method.  The problem easily reduces to proving uniform $C^{\infty}$ \emph{a priori} estimates for a smooth solution $\varphi$ of (\ref{ma0}), depending only on $(M,g)$ and $F$.  It is natural to ask:  to what extent can Yau's results be generalized to the case when $\omega$ is no longer K\"ahler?  

The current paper was motivated by the recent  interesting work of B. Guan and Q. Li.   In \cite{GL},  Guan-Li made a number of  important advances, including the result that 
\begin{equation} \label{gl}
\partial \ov{\partial} \omega^k =0, \quad k=1,2,
\end{equation} 
is a sufficient condition for solving equation (\ref{ma0}).  Note that (\ref{gl}) implies that $\partial \ov{\partial} \omega^k =0$ for all $1 \le k \le n-1$ (see, for example, \cite{FT}).
Also, Guan-Li apply their estimates to the problem of finding geodesics in the space of Hermitian metrics via the homogeneous complex Monge-Amp\`ere equation   (for some related works, see  \cite{CLN, BT, Ma, Se, D1, GuB, Ch, GuP, PS}, for example).

We focus now on the question of obtaining \emph{a priori} estimates for solutions $\varphi$ of (\ref{ma0}).  In complex dimension two, we show that \emph{no assumption} on $(M, \omega)$ is needed to obtain the \emph{a priori} estimates.  In higher dimensions we impose the  condition
\begin{equation} \label{balance}
d (\omega^{n-1})=0.
\end{equation}
Hermitian metrics $\omega$ that satisfy this equality are known as balanced metrics.  The condition (\ref{balance}) appears to be a natural one  and  has been much studied, due in part to a role in string theory.   We refer the reader to \cite{AB, BBFTY, FIUV, FY, FLY, Mi} and the references therein for more details.  

We give a proof of the following:

\bigskip
\noindent
{\bf Main Theorem} \, \emph{Let $(M,\omega)$ be a compact Hermitian manifold of complex dimension $n$. Let $\varphi$ be a smooth solution of the complex Monge-Amp\`ere equation (\ref{ma0}), subject to (\ref{sup}).  Assume that either
\begin{enumerate}
\item[(a)] $n=2$; \ \ or
\item[(b)] $n>2$ and $d(\omega^{n-1})=0$.
\end{enumerate}
Then there are uniform $C^{\infty}$ a priori estimates on $\varphi$ depending only on $(M,\omega)$ and  $F$.}

\bigskip

Shortly after writing the first version of this paper, the authors discovered that the result of the Main Theorem is contained in a paper of Cherrier \cite{ChP}.  However, since our methods differ from those of Cherrier  at almost every stage of the argument, we believe that our alternative proof of the Main Theorem may still be of interest.  We describe later in the introduction some new elements of our proof.

First we make a couple of small remarks.  We note that it is not necessary to assume that the solution is \emph{a priori} smooth ($\varphi\in C^{2}$ would be enough, for example), but we have made this assumption for the sake of simplicity.
Also, with some work, one could make explicit the dependence of the estimates on $\varphi$ in terms of $(M,g)$ and the bounds $\| F \|_{C^k}$ for $k \in \mathbb{N}$.  

\begin{remark}  Very recently, B. Guan and Q. Li have informed us that they can also prove case (a) of the Main Theorem by extending their methods. 
\end{remark}

We now give two corollaries of the Main Theorem.  The first is contained in \cite{ChP} and states that one can solve (\ref{ma0}) after possibly adding a constant to the function $F$.

\bigskip
\noindent
{\bf Corollary 1} \, \emph{Let $(M, \omega)$ be a compact Hermitian manifold satisfying one of the conditions (a) or (b) of the Main Theorem.  Then for any  smooth function $F$ on $M$ there exists a constant $b$ and a smooth function $\varphi$ on $M$ solving
\begin{equation} \label{ma3}
(\omega + \sqrt{-1} \partial \ov{\partial} \varphi)^n = e^{F+b} \omega^n.
\end{equation}}

If $\varphi$ satisfies the normalization condition (\ref{sup}) then the $\varphi$ and $b$ are unique (see Remark \ref{unique}.1).
Cherrier's proof of Corollary 1 makes use of a rather general result of Delano\"e \cite{De}.  
For the sake of completeness, we include a direct proof of Corollary 1 in Section \ref{sectionopenness}. 
  
Note that  the result of Guan-Li   does not require the addition of a constant $b$ to the function $F$. The basic reason for this difference is that, unlike the balanced condition, the assumption (\ref{gl}) implies
$$\int_M (\omega + \sqrt{-1} \partial \ov{\partial} \varphi)^n = \int_M \omega^n.$$
We also remark that, in general, the balanced condition is in some sense disjoint from the condition (\ref{gl}) of Guan-Li.  Indeed, if a balanced metric $\omega$ satisfies $\de\db\omega=0$ it has to be K\"ahler (see e.g. \cite[Proposition 1.4]{FPS}). There are many known examples of non-K\"ahler balanced metrics \cite{AB, FIUV, FLY, Mi}, that hence do not satisfy \eqref{gl}.   There seem to be fewer examples in the literature of non-K\"ahler metrics satisfying \eqref{gl} in dimensions $n>2$.  However,  an elementary such example is the product of a complex curve with a K\"ahler metric and a complex surface with a non-K\"ahler Gauduchon metric. More examples were recently constructed in \cite{FT}.  

We also point out  that Cherrier proves the Main Theorem (and thus Corollary 1) under a technical assumption on $\omega$  which is slightly weaker than balanced and can be regarded as `close to balanced'.   He also considers the case of conformally K\"ahler $\omega$.

Next we give an application of Corollary 1 to the first Chern class of the Hermitian manifold $(M,\omega)$.  Recall that for K\"ahler manifolds, Yau's theorem implies the Calabi Conjecture, which states that any representative of the first Chern class can be written as the Ricci curvature of a K\"ahler metric in any given K\"ahler class.  We restrict now to two complex dimensions, and show that 
 an analogue of the Calabi Conjecture is true  if and only if an integral condition  holds.  Denoting by $\textrm{Ric}(\omega)$ the first Chern form of the canonical connection of $g$ (see Section \ref{sectionchern}), we show that:

\pagebreak[3]
\bigskip
\noindent
{\bf Corollary 2} \,  \emph{Let $(M, \omega)$ be a compact Hermitian manifold of complex dimension two. Denote by $\omega_G$ any Gauduchon metric on $M$ (see (\ref{gauduchon})). 
Given a closed real $(1,1)$ form $\psi$ cohomologous to $c_1(M)$, 
we can write 
\begin{equation}\label{eqn}
\psi=\Ric(\omega+\sqrt{-1}\de\db\varphi),
\end{equation}
for some Hermitian metric $\omega+\sqrt{-1}\de\db\varphi$ if and only if
\begin{equation}\label{constr}
\int_M (\Ric(\omega)-\psi)\wedge\omega_G=0.
\end{equation}}

We now briefly outline the strategy for our proof of the Main Theorem and point out how our approach differs from that of \cite{ChP} and also \cite{GL}.
We begin by describing the methods of \cite{ChP}, \cite{GL}.   Both  papers give an estimate on the quantity 
$\tr{g}{g'}$ in terms of $\varphi$, where $g'$ is the Hermitian metric defined  by
\begin{equation}\label{herm}
g'_{i \ov{j}} = g_{i \ov{j}} + \partial_i \partial_{\ov{j}} \varphi>0.
\end{equation}
The precise form of the estimate on $\tr{g}{g'}$ is not important for the methods of \cite{ChP} or \cite{GL}, but
tracing through the arguments of \cite{GL} one can see for example that
\begin{equation} \label{chestimate}
\tr{g}{g'} \le Ce^{\left(e^{A(\sup_M\varphi- \inf_M \varphi)} - e^{A(\sup_M\varphi- \varphi)}\right)}.
\end{equation}
Note that (\ref{chestimate}) is proved without any conditions on $(M, \omega)$.   Zhang \cite{Zh} independently proved a similar estimate on $\tr{g}{g'}$ and he also considered a generalization of the complex Monge-Amp\`ere equation.   

From (\ref{chestimate}), the metric $g'$ is uniformly bounded once we have a bound on $\varphi$.  As remarked in \cite{GL}, it is not clear whether the $C^0$ estimate of $\varphi$ given in \cite{Y}, or the alternative methods of \cite{Ko} or \cite{Bl1}, can be extended to the general Hermitian case.  But by making use of the assumptions (a) or (b) and (\ref{gl}) respectively, a $C^0$ bound on $\varphi$ is proved in \cite{ChP} and \cite{GL}.  The argument for the $C^0$ estimate in \cite{ChP} is particularly intricate, especially in the case  (b).

Our approach is quite different.  We also prove an estimate on $\tr{g}{g'}$ depending on $\varphi$, but we show that the assumptions (a) and (b) imply an  improved  version of (\ref{chestimate}) and then make use of this to   obtain the $C^0$  bound of $\varphi$.    We prove the following estimate in Section \ref{sectionso}.

\bigskip
\noindent
{\bf Theorem 2.1} \ \emph{Under the assumptions of the Main Theorem, there exist uniform constants $C$ and $A$ such that 
\begin{equation} \label{so2}
\emph{tr}_{g}{g'} \le C e^{A (\varphi - \inf_M \varphi)}.
\end{equation}}

We remark that our proof of this improved estimate was motivated by our reading the paper \cite{GL}.  We  use and modify some key ideas from \cite{GL} (variations of some of these ideas also appear in \cite{ChP}, \cite{Zh}).

Observe that (\ref{so2}) is now of the same form as the second order estimate from \cite{Y}. 
We then prove the uniform $C^0$ bound on $\varphi$ by
  making use of an idea from \cite{W1, W2} (see also \cite{TWY}), where it was shown that in the K\"ahler case, the estimate  (\ref{so2}) implies a uniform bound on $\varphi$.
Indeed, we show in Section \ref{sectionzo} that combining (\ref{so2}) with a Moser iteration argument applied to the exponential of $\varphi$ gives,
\begin{equation} \label{a}
- \inf \vp \le C_{\alpha}+ \log \left( \frac{1}{\int_M \omega^n} \int_M e^{-\alpha \vp} \omega^n \right)^{\frac{1}{\alpha}},
\end{equation}
for all $\alpha>0$.   Some further work shows that (\ref{a}) gives  a uniform bound of $\| \varphi \|_{C^0}$ and hence $\tr{g}{g'}$.  We note that our $C^0$ estimate follows from (\ref{so2}) and does not make use of the complex Monge-Amp\`ere equation.  It would be interesting to know whether one can prove an estimate similar to (\ref{so2}) for other equations on Hermitian manifolds, and then apply our $C^0$ estimate.

\pagebreak[3]
\begin{remark} A consequence of our results is that if the improved second order estimate (\ref{so2}) could be shown to hold in the case of general Hermitian $\omega$, then one could remove the conditions (a), (b) in the statements of the Main Theorem and Corollary 1.  Indeed, our argument for the zero order  estimate can easily be seen to hold without assuming (a) or (b).  The higher order estimates also do not require these assumptions.  
\end{remark}

In Section \ref{sectionproof} we complete the proof of the Main Theorem.  
Once we have the bound on $\tr{g}{g'}$ it is a straightforward matter to finish the proof. Cherrier \cite{ChP} obtains higher order estimates via  a third order estimate on $\varphi$ analogous to that of \cite{Y}, \cite{Ca2}.  It was shown in \cite{GL} that the higher order estimates follow from a bound on the real Hessian of $\varphi$ and arguments of Evans  \cite{Ev} and Krylov \cite{Kr}.   In order to give a short and largely self-contained proof we adapt instead an 
 estimate of Trudinger \cite{Tr}.  
 
In Sections \ref{sectionopenness} and \ref{sectionchern}  we prove Corollaries 1 and 2 respectively.

\begin{remark}  
It is perhaps worth pointing out that a different generalization of Yau's theorem to non-K\"ahler manifolds has recently been proposed by Donaldson \cite{D2}.  If $\omega$ is instead a symplectic form on a compact real $4$-dimensional manifold, compatible with an almost complex structure $J$, then one can ask whether there exists a symplectic form $\tilde{\omega}$ compatible with $J$ solving
$$\tilde{\omega}^2 = e^F \omega^2.$$
This problem has a rather different nature to the study of (\ref{ma0}), since in general one does not expect that $\tilde{\omega}$ is of the form $\omega + d(Jd \varphi)$.  For some recent developments on this, we refer the reader to \cite{W2, TWY, TW1, TW2}.

\end{remark}

\section{The second order estimate} \label{sectionso}
\setcounter{equation}{0}

Write $g'$ for the Hermitian metric given by (\ref{herm}).
In this section we prove the following estimate on $\tr{g}{g'} = g^{i \ov{j}}g'_{i \ov{j}}$, depending on $\varphi$.

\begin{theorem} \label{tso}
Under the assumptions of the Main Theorem, there exist uniform constants $C$ and $A$ such that 
\begin{equation} \label{so}
\emph{tr}_{g}{g'} \le C e^{A (\varphi - \inf_M \varphi)}.
\end{equation}
\end{theorem}

Note that here, and henceforth, when we say `uniform constant' it should be understood that we mean a constant depending only on $(M,g)$ and bounds for $F$.   We will often write $C$ or $C'$ for such a constant, where the value of $C$, $C'$ may differ from line to line.

Write $\Delta$ for the Laplace operator of the canonical connection (i.e. the Chern connection) associated to the Hermitian metric $g$.   Namely, for a function $f$,
$$\Delta f = g^{i \ov{j}} \partial_i \partial_{\ov{j}} f=\frac{n\omega^{n-1}\wedge\sqrt{-1}\de\db f}{\omega^n},$$
and similarly for $\Delta'$.  Note that $\tr{g}{g'} = n+ \Delta \vp$ and $\tr{g'}{g} = n - \Delta' \varphi$.

\begin{proof}[Proof of Theorem \ref{tso}]  
Following the basic method of Yau \cite{Y} (see also \cite{A}),  we apply the maximum principle to the quantity
$$Q = \log \tr{g}{g'} -A \varphi,$$
for some suitably chosen large constant $A$.  For convenience, we will write $E_1$ and $E_2$ for terms satisfying
$$|E_1| \le C \tr{g'}{g}, \quad |E_2| \le C' (\tr{g}{g'})( \tr{g'}{g}).$$
Note that the equation (\ref{ma0}), which can be written as, $$\det{g'} = e^F \det{g},$$ implies that $\tr{g}{g'}$ and $\tr{g'}{g}$ are uniformly bounded from below away from zero.  It follows that  a uniform constant is itself of type $E_1$  and any term of type $E_1$ is also of type $E_2$. 
We also note here the elementary inequality
\begin{equation}\label{triv}
\tr{g}{g'}\leq \frac{1}{(n-1)!} (\tr{g'}{g})^{n-1}\frac{\det{g'}}{\det{g}}= \frac{1}{(n-1)!} (\tr{g'}{g})^{n-1}e^F,
\end{equation}
and in the case of dimension $n=2$, the equality
\begin{equation} \label{triv2}
\tr{g}{g'} = (\tr{g'}{g})e^F.
\end{equation}

We will show that in both cases (a) and (b) from the statement of the Main Theorem, we have the estimate
\begin{equation} \label{key}
\Delta' \log \tr{g}{g'} \ge E_1.
\end{equation}
Given (\ref{key}), the theorem follows immediately.  Indeed, if $Q$ achieves its maximum at a point $p$, we obtain
$$0 \ge (\Delta' Q) (p) \ge E_1 - An + A(\tr{g'}{g})(p) \ge (\tr{g'}{g})(p) -C,$$
if $A$ is chosen sufficiently large.  Then $(\tr{g'}{g})(p)$ is uniformly bounded from above, which implies by \eqref{triv} that $(\tr{g}{g'})(p)$ is too.  Then for any $q \in M$,
$$Q(q) \le Q(p) \le \log C - A\inf_M \varphi,$$
and (\ref{so}) follows after exponentiating.

We need the following lemma from \cite{GL} (see also \cite{ST} for a similar argument).

\begin{lemma} \label{lemmagl}  At any point $p$ there exists a holomorphic coordinate system centered at $p$ such that for all $i,j$, 
\begin{equation} \label{condition}
g_{i \ov{j}}(0) = \delta_{ij}, \quad \partial_{j} g_{i \ov{i}}(0) =0,
\end{equation}
and also that  the matrix $(\de_k \de_{\ov{\ell}}\vp)(0)$ is diagonal.
\end{lemma}

When we perform local computations we will always assume we are using a coordinate system given by this lemma at a point $p$, and we will use lower indices to denote partial derivatives.
We first compute:
\begin{eqnarray} \nonumber
\Delta' \tr{g}{g'} & = & g'^{i \ov{j}} \dd{i} \dd{\ov{j}} (g^{k\ov{\ell}} \vp_{k\ov{\ell}}) \\ \nonumber
& = & \sum_{i,k} g'^{i \ov{i}} \vp_{k \ov{k} i \ov{i}} - 2 \textrm{Re} \left( \sum_{i,j,k} g'^{i \ov{i}} \partial_{\ov{i}} g_{j \ov{k}} \vp_{k \ov{j} i} \right) + \sum_{i, j,k} g'^{i \ov{i}} \partial_i g_{j \ov{k}} \partial_{\ov{i}} g_{k \ov{j}} \vp_{k \ov{k}} \\ \nonumber
&&\mbox{} + \sum_{i,j,k} g'^{i \ov{i}} \partial_{i} g_{k \ov{j}} \partial_{\ov{i}} g_{j \ov{k}} \vp_{k \ov{k}} - \sum_{i,k} g'^{i\ov{i}} \partial_{i} \partial_{\ov{i}} g_{k\ov{k}} \vp_{k \ov{k}} \\ \label{deltatr}
& = &  \sum_{i,k} g'^{i \ov{i}} \vp_{k \ov{k} i \ov{i}} - 2 \textrm{Re} \left( \sum_{i,j,k} g'^{i \ov{i}} \partial_{\ov{i}} g_{j \ov{k}} \vp_{k \ov{j} i} \right)  + E_2.
\end{eqnarray}
We now apply the operator $\Delta \log (\cdot)$ to the equation (\ref{ma0}), to get:
\begin{eqnarray} \nonumber
\lefteqn{- g^{k \ov{\ell}} g'^{p \ov{j}} g'^{i\ov{q}} \partial_k g'_{p \ov{q}} \partial_{\ov{\ell}} g'_{i \ov{j}} + g^{k \ov{\ell}} g'^{i \ov{j}} \partial_k \partial_{\ov{\ell}} g'_{i \ov{j}} } \\& = & \Delta F - g^{k\ov{\ell}} g^{p \ov{j}} g^{i \ov{q}} \partial_k g_{p\ov{q}} \partial_{\ov{\ell}} g_{i \ov{j}} + g^{k \ov{\ell}} g^{i \ov{j}} \partial_k \partial_{\ov{\ell}} g_{i \ov{j}}.
\end{eqnarray}
We can rewrite this as:
\begin{equation} \label{eq1}
\sum_{i,k} g'^{i \ov{i}} \vp_{i \ov{i} k \ov{k}} = \sum_{i,j,k} g'^{i \ov{i}} g'^{j \ov{j}} \partial_{k} g'_{i \ov{j}} \partial_{\ov{k}} g'_{j \ov{i}} + E_1.
\end{equation}
Thus (\ref{deltatr}) and (\ref{eq1}) give
\begin{equation} \label{deltatr2}
\Delta' \tr{g}{g'} =  \sum_{i,j,k} g'^{i \ov{i}} g'^{j \ov{j}} \partial_{k} g'_{i \ov{j}} \partial_{\ov{k}} g'_{j \ov{i}} - 2 \textrm{Re} \left( \sum_{i,j,k} g'^{i \ov{i}} \partial_{\ov{i}} g_{j \ov{k}} \vp_{k \ov{j} i} \right)  + E_2.
\end{equation}

A key trick in \cite{GL}  is to use part of the first term $\sum_{i,j,k} g'^{i \ov{i}} g'^{j \ov{j}} \partial_{k} g'_{i \ov{j}} \partial_{\ov{k}} g'_{j \ov{i}}$ on the right hand side of (\ref{deltatr2}) to deal with the troublesome second term.  Calculate, using (\ref{condition}),
\begin{eqnarray} \nonumber
\sum_{i,j,k} g'^{i\ov{i}} \partial_{\ov{i}} g_{j \ov{k}} \vp_{k \ov{j} i} & = & \sum_{i,j,k} g'^{i \ov{i}} \partial_{\ov{i}} g_{j \ov{k}} ( \partial_k g_{i \ov{j}} + \vp_{i \ov{j} k}) + E_1 \\
&= & \sum_{i} \sum_{j \neq k} g'^{i\ov{i}} \partial_{\ov{i}} g_{j \ov{k}} \partial_{k} g'_{i \ov{j}} + E_1.
\end{eqnarray}
It follows that
\begin{eqnarray} \nonumber
\lefteqn{ \left| 2 \textrm{Re} \left( \sum_{i,j,k} g'^{i \ov{i}} \partial_{\ov{i}} g_{j \ov{k}} \vp_{k \ov{j} i} \right) \right|  } \\ \nonumber & \le & \frac{1}{2} \sum_i \sum_{j \neq k} g'^{i \ov{i}} g'^{j \ov{j}} \partial_k g'_{i \ov{j}} \partial_{\ov{k}} g'_{j \ov{i}} + 2\sum_i \sum_{j \neq k} g'^{i \ov{i}} g'_{j \ov{j}} \partial_{\ov{i}} g_{j \ov{k}} \partial_i g_{k \ov{j}} +E_1 \\ \label{eq2}
& = & \frac{1}{2} \sum_i \sum_{j \neq k} g'^{i \ov{i}} g'^{j \ov{j}} \partial_k g'_{i \ov{j}} \partial_{\ov{k}} g'_{j \ov{i}}  + E_2.
\end{eqnarray}
Note that the factor $1/2$ in (\ref{eq2}) differs from the corresponding factor $1$ in \cite{GL}.  We will make use of this later.
Combining (\ref{deltatr2}) and (\ref{eq2}) we have
\begin{eqnarray} \label{dt}
\Delta' \tr{g}{g'} \ge \sum_{i,j} g'^{i\ov{i}} g'^{j \ov{j}} \partial_j g'_{i \ov{j}} \partial_{\ov{j}} g'_{j \ov{i}} +  \frac{1}{2} \sum_i \sum_{j \neq k} g'^{i \ov{i}} g'^{j \ov{j}} \partial_k g'_{i \ov{j}} \partial_{\ov{k}} g'_{j \ov{i}}  + E_2.
\end{eqnarray}

We now claim that
\begin{equation} \label{hope}
\frac{| \partial \tr{g}{g'}|^2_{g'}}{\tr{g}{g'}} \le  \sum_{i,j} g'^{i \ov{i}} g'^{j\ov{j}} \partial_j g'_{i \ov{j}} \partial_{\ov{j}} g'_{j \ov{i}}+ \frac{1}{2} \sum_i \sum_{j \neq k} g'^{i \ov{i}} g'^{j \ov{j}} \partial_k g'_{i \ov{j}} \partial_{\ov{k}} g'_{j \ov{i}} + E_2.
\end{equation}
Given (\ref{hope}) we would then have (\ref{key}) since
$$\Delta' \log \tr{g}{g'} = \frac{\Delta' \tr{g}{g'}}{\tr{g}{g'}} - \frac{| \partial \tr{g}{g'}|^2_{g'}}{(\tr{g}{g'})^2} \ge E_1.$$  
Thus the proof of the theorem reduces to establishing (\ref{hope}).  We split the proof of (\ref{hope}) into the two cases (a) and (b).  

We remark that a different inequality from (\ref{hope}) is proved in \cite{GL}, where they make use of an alternative expression for $Q$.   Indeed, it is not clear whether (\ref{hope}) holds in general without our assumptions (a) or (b).

\bigskip
\noindent
{\it Proof of (\ref{hope}) in the case (a).} \, 
We prove (\ref{hope}) assuming that $n=2$.  First, making use of (\ref{condition}), we have
\begin{equation} \label{dtr}
\partial_i \tr{g}{g'} = \sum_{j} \partial_i g'_{j \ov{j}}.
\end{equation}
We use an argument in Yau's second order estimate \cite{Y}.   Applying the  Cauchy-Schwarz inequality twice:
\begin{eqnarray} \nonumber
\frac{| \partial \tr{g}{g'} |^2_{g'}}{\tr{g}{g'}}& = & \frac{1}{\tr{g}{g'}} \sum_{i,j,k} g'^{i \ov{i}} \partial_i g'_{j \ov{j}} \partial_{\ov{i}} g'_{k \ov{k}} \\ \nonumber
& = & \frac{1}{\tr{g}{g'}} \sum_{j,k} \sum_i \sqrt{g'^{i\ov{i}}}  \partial_i g'_{j \ov{j}} \sqrt{g'^{i\ov{i}}} \partial_{\ov{i}} g'_{k \ov{k}} \\ \nonumber
& \le & \frac{1}{\tr{g}{g'}} \sum_{j,k} \left( \sum_i g'^{i \ov{i}} | \partial_i g'_{j \ov{j}}|^2 \right)^{1/2} \left( \sum_i g'^{i \ov{i}} |\partial_i g'_{k \ov{k}}|^2 \right)^{1/2} \\ \nonumber
& = & \frac{1}{\tr{g}{g'}} \left( \sum_j \left( \sum_i g'^{i \ov{i}} | \partial_i g'_{j \ov{j}}|^2 \right)^{1/2} \right)^2 \\ \nonumber
& = & \frac{1}{\tr{g}{g'}} \left( \sum_j \sqrt{g'_{j \ov{j}}} \left( \sum_i g'^{i \ov{i}} g'^{j \ov{j}} | \partial_i g'_{j \ov{j}}|^2 \right)^{1/2} \right)^2 \\
& \le & \sum_{i,j} g'^{i \ov{i}} g'^{j \ov{j}}  \partial_i g'_{j \ov{j}} \partial_{\ov{i}} g'_{j \ov{j}}. \label{siu}
\end{eqnarray}

Using (\ref{condition}) again, we can rewrite this as,
\begin{eqnarray} \nonumber
\frac{| \partial \tr{g}{g'} |^2_{g'}}{\tr{g}{g'}} & \le & \sum_{i,j} g'^{i \ov{i}} g'^{j \ov{j}} \vp_{j \ov{j} i} \vp_{j \ov{j} \ov{i}}  \\ \nonumber
& = & \sum_{i,j} g'^{i \ov{i}} g'^{j \ov{j}} \vp_{i \ov{j} j} \vp_{j \ov{i} \ov{j}} \\ \nonumber
& = &  \sum_{i,j} g'^{i \ov{i}} g'^{j \ov{j}} (\partial_j g'_{i \ov{j}} - \partial_j g_{i \ov{j}})( \partial_{\ov{j}} g'_{j \ov{i}} - \partial_{\ov{j}} g_{j \ov{i}} ) \\ \nonumber
& = & \sum_{i,j} g'^{i \ov{i}} g'^{j \ov{j}} \partial_j g'_{i \ov{j}}  \partial_{\ov{j}} g'_{j \ov{i}} - 2\textrm{Re} \left(\sum_{i,j} g'^{i \ov{i}} g'^{j \ov{j}} \partial_j g'_{i \ov{j}}  \partial_{\ov{j}} g_{j \ov{i}}  \right) \\ \nonumber
&&\mbox{}+ \sum_{i,j} g'^{i \ov{i}} g'^{j \ov{j}} \partial_j g_{i \ov{j}}\partial_{\ov{j}} g_{j \ov{i}}\\ \nonumber
& = & \sum_{i,j} g'^{i \ov{i}} g'^{j \ov{j}} \partial_j g'_{i \ov{j}}  \partial_{\ov{j}} g'_{j \ov{i}} - 2\textrm{Re} \left(\sum_{i,j} g'^{i \ov{i}} g'^{j \ov{j}} \partial_i g'_{j \ov{j}}  \partial_{\ov{j}} g_{j \ov{i}}  \right) \\ \nonumber
&&\mbox{}- \sum_{i,j} g'^{i \ov{i}} g'^{j \ov{j}} \partial_j g_{i \ov{j}}\partial_{\ov{j}} g_{j \ov{i}} \\ 
& \leq & \sum_{i,j} g'^{i \ov{i}} g'^{j \ov{j}} \partial_j g'_{i \ov{j}}  \partial_{\ov{j}} g'_{j \ov{i}} - 2\textrm{Re} \left(\sum_{i \neq j} g'^{i \ov{i}} g'^{j \ov{j}} \partial_i g'_{j \ov{j}}  \partial_{\ov{j}} g_{j \ov{i}}  \right).
\label{s1}
   \end{eqnarray}
It remains to control the term
$$2\textrm{Re} \left(\sum_{i \neq j} g'^{i \ov{i}} g'^{j \ov{j}} \partial_i g'_{j \ov{j}}  \partial_{\ov{j}} g_{j \ov{i}}  \right).$$
Using again (\ref{condition}) and also the fact that in dimension $n=2$, $\tr{g}{g'}$ and $\tr{g'}{g}$ are uniformly equivalent, which follows from \eqref{triv2}, we see that
\begin{eqnarray} \nonumber
\left| 2\textrm{Re} \left(\sum_{i \neq j} g'^{i \ov{i}} g'^{j \ov{j}} \partial_i g'_{j \ov{j}}  \partial_{\ov{j}} g_{j \ov{i}}  \right) \right| & \le & \frac{1}{2} \sum_{i \neq j} g'^{j \ov{j}} g'^{j\ov{j}} | \partial_i g'_{j \ov{j}}|^2+ 2 \sum_{i \neq j} g'^{i \ov{i}}g'^{i \ov{i}} | \partial_j g_{i\ov{j}}|^2   \\  \label{s2}
& = &\frac{1}{2} \sum_{i \neq j} g'^{j \ov{j}} g'^{j\ov{j}} | \partial_i g'_{j \ov{j}}|^2  + E_2.
\end{eqnarray}
But notice that 
\begin{equation} \label{s3}
\frac{1}{2}\sum_{i \neq j} g'^{j \ov{j}} g'^{j\ov{j}} | \partial_i g'_{j \ov{j}}|^2 \le \frac{1}{2} \sum_i \sum_{j \neq k} g'^{i \ov{i}} g'^{j \ov{j}} \partial_k g'_{i \ov{j}} \partial_{\ov{k}} g'_{j \ov{i}}.
\end{equation}
Combining (\ref{s1}), (\ref{s2}) and (\ref{s3}) gives (\ref{hope}).

\bigskip
\noindent
{\it Proof of (\ref{hope}) in the case (b).} \, We finish the proof of the theorem by establishing (\ref{hope}) when $n>2$ and $d(\omega^{n-1})=0$.  The condition $d(\omega^{n-1})=0$ is equivalent to $\sum_j T^j_{ji}=0$, where $T^k_{ji}=g^{k\ov{\ell}}(\de_j g_{i\ov{\ell}}-\de_i g_{j\ov{\ell}})$ is the torsion of $g$ (see e.g. \cite{G2, Ga2}).  Then we have, at $p$, using the coordinate system of Lemma \ref{lemmagl},
\begin{equation} \label{balanced}
0= \sum_j T^j_{ji}= \sum_j \partial_j g_{i\ov{j}} - \sum_j \partial_i g_{j \ov{j}}= \sum_j \partial_j g_{i \ov{j}}.
\end{equation}
Then from (\ref{dtr}) we have,
\begin{equation}
\partial_i \tr{g}{g'} = \sum_j \partial_i g'_{j \ov{j}} = \sum_j \partial_i \varphi_{j\ov{j}} = \sum_j (\partial_j g_{i\ov{j}} + \partial_j \varphi_{i \ov{j}}) = \sum_j \partial_j g'_{i \ov{j}}.
\end{equation}
Repeating the argument we used in (\ref{siu}) we  obtain
\begin{eqnarray} \nonumber
\frac{| \partial \tr{g}{g'} |^2_{g'}}{\tr{g}{g'}}& = & \frac{1}{\tr{g}{g'}} \sum_{i,j,k} g'^{i \ov{i}} \partial_j g'_{i \ov{j}} \partial_{\ov{k}} g'_{k \ov{i}} \\ 
& \le & \sum_{i,j} g'^{i \ov{i}} g'^{j \ov{j}}  \partial_j g'_{i \ov{j}} \partial_{\ov{j}} g'_{j \ov{i}}, \label{siu2}
\end{eqnarray}
and this immediately gives (\ref{hope}).  Q.E.D.
\end{proof}

\section{The zeroth order estimate} \label{sectionzo}
\setcounter{equation}{0}

In this section we will prove:

\begin{theorem} \label{tzo}
Under the assumptions of the Main Theorem,  there exists a uniform constant $C$  such that 
\begin{equation} \label{zo}
\| \varphi \|_{C^0} \le C.\end{equation}
\end{theorem}

\begin{proof}  We note that the proof of Theorem \ref{tzo} crucially requires the precise form of Theorem \ref{tso}, but  does not again make use of the equation (\ref{ma0}). In particular, if one could prove \eqref{so} without assumptions (a) or (b), Theorem \ref{tzo} would also follow.
 We first prove two lemmas which we will use in both cases (a) and (b). 
 \begin{lemma} \label{lemmaalpha}
For every $\alpha>0$ there exists a constant $C_{\alpha}$ such that
\begin{equation} \label{alpha}
- \inf_M \vp \le C_{\alpha}+ \log \left( \int_M e^{-\alpha \vp} d\mu \right)^{\frac{1}{\alpha}},
\end{equation}
where $d\mu = \omega^n/\int_M \omega^n$.
\end{lemma}
\begin{proof}  The argument is almost identical  to the one given in \cite{W2}.
Let $\delta>0$ be a given small constant.  
Define $w = e^{-B \vp}$, $\gamma = 1- \delta$, $B=A/\gamma$, where $A$ is the constant given in (\ref{so}).  For $p\ge 1$, compute:
\begin{eqnarray*}
\int_M | \partial w^{p/2}|^2_g \omega^n & = & n \int_M  \sqrt{-1} \partial e^{-\frac{Bp\vp}{2}} \wedge \ov{\partial} e^{-\frac{Bp\vp}{2}} \wedge \omega^{n-1} \\
& \le & C n p^2  \int_M e^{-Bp\vp} \sqrt{-1} \partial \vp \wedge \ov{\partial} \vp  \wedge  \omega^{n-1} \\
& \le & - C n p \int_M \sqrt{-1} \partial \left( e^{-Bp \vp }\right) \wedge \ov{\partial} \vp \wedge \omega^{n-1} \\
& = & C np  \int_M e^{-Bp \vp } \sqrt{-1} \partial \ov{\partial} \vp \wedge \omega^{n-1} \\ 
&& \mbox{} - C np  \int_M e^{-Bp\vp } \sqrt{-1} \, \ov{\partial} \vp \wedge \partial \omega^{n-1} \\
& = & C p  \int_M e^{-Bp \vp } \Delta \vp \, \omega^n \\ 
&& \mbox{} + C'  \int_M \sqrt{-1} \, \ov{\partial} \left( e^{-Bp \vp } \right)  \wedge \partial \omega^{n-1} \\
& \le & C p \| w \|^{\gamma}_{C^0} \int_M e^{-B(p- \gamma)\vp} \omega^n \\
&& \mbox{} + C p \int_M \sqrt{-1} e^{-B(p- \gamma) \vp} e^{-B\gamma \vp}  \partial \ov{\partial} (\omega^{n-1}) \\
& \le & C p \| w \|^{\gamma}_{C^0} \int_M e^{-B(p-\gamma) \vp} \omega^n\\
& = & C p \| w \|^{\gamma}_{C^0} \int_M w^{p-\gamma} \omega^n,
\end{eqnarray*}
where we have used the fact that $\partial \ov{\partial} (\omega^{n-1})$ is bounded and
$$\Delta \vp=\tr{g}{g'}-n\le C e^{A(\vp-\inf \vp)} = C e^{B \gamma (\vp - \inf \vp)} \le C \| w \|^{\gamma}_{C^0} e^{B \gamma \vp}.$$

Now, the Sobolev inequality gives us, for $\beta =n/(n-1)$, and any smooth $f$,
\begin{equation}\label{sob}
\left( \int_M f^{2\beta} \omega^n \right)^{1/\beta} \le C \left( \int_M | \partial f |_g^2 \omega^n + \int_M f^2 \omega^n \right),
\end{equation}
%
%
%
Applying this to $f=w^{p/2}$ and raising to the power $1/p$, we obtain
\begin{equation} \label{eqniterate}
\| w \|_{p\beta} \le C^{1/p} p^{1/p} \| w \|_{C^0}^{\gamma/p} \| w \|_{p-\ga}^{(p-\ga)/p},
\end{equation}
where $\| \ \|_q$ denotes the $L^q$ norm with respect to $d\mu$ (allowing $0<q<1$, defined in the obvious way).
By the same iteration as in \cite{W2}  (see also \cite{TWY}) we replace $p$ with $p\beta + \gamma$ in (\ref{eqniterate}) to obtain for $k=1,2, \ldots$,
\begin{equation} \label{eqnw}
 \| w \|_{p_k \beta} \le C \| w \|_{C^0}^{1-a(k)} \| w \|_{p-\gamma}^{a(k)},
 \end{equation}
where $p(k) \rightarrow \infty$ and $a(k) \rightarrow a \in (0,1)$ as $k \rightarrow \infty$.  We refer the reader to \cite{W2} for the details. Letting $k \rightarrow \infty$ and setting $p=1$ in (\ref{eqnw}), we have
$$\| w \|_{C^0} \le C \| w \|_{\delta}.$$
Choosing $\delta$ sufficiently small proves the lemma for $\alpha$ sufficiently small.  It follows from H\"older's inequality  that if Lemma \ref{lemmaalpha} holds for all $\alpha>0$ sufficiently small then it holds for all $\alpha>0$. Q.E.D.

\end{proof}

Denote by $|\cdot|$ the measure of a set with respect to the measure $d\mu = \omega^n/\int_M \omega^n$.   Then we have the following.

\begin{lemma} \label{lemmameasure}
Let $C_1$ be the constant corresponding to $\alpha=1$ from Lemma \ref{lemmaalpha}.  Then
\begin{equation} \label{enough}
|\{\varphi\leq \inf_M \varphi +C_1+1\}|\geq\frac{e^{-C_1}}{4}.
\end{equation}
\end{lemma}

\begin{proof}  Suppose for a contradiction that 
\begin{equation} \label{ee3}
|\{\varphi\leq \inf_M \varphi +C_1+1\}|<\frac{e^{-C_1}}{4}.
\end{equation}
First note that from Lemma \ref{lemmaalpha} with $\alpha=1$, we have
\begin{equation} \label{ee1}
e^{-\varphi} \leq e^{-\inf_M\varphi}\leq e^{C_1} \int_M e^{-\varphi},
\end{equation}
where here and henceforth we are integrating with respect to $d\mu$.
Now
\begin{equation} \label{ee2}
\int_M e^{-\varphi} = \int_{ \{ e^{-\varphi} \ge \frac{1}{e} \int_M e^{-\varphi} \}} e^{-\varphi}  + \int_{ \{e^{-\varphi} < \frac{1}{e} \int_M e^{-\varphi} \} }e^{-\varphi}.
\end{equation}
But
\begin{eqnarray*}
| \{ e^{-\varphi} \ge \frac{1}{e} \int_M e^{-\varphi} \} | & \le & | \{ e^{-\varphi} \ge\frac{1}{e^{C_1+1}} e^{-\inf_M\varphi} \} | \\
& = & | \{ - \varphi \ge -\inf_M\varphi-C_1-1 \} | \\
& < & \frac{e^{-C_1}}{4},
\end{eqnarray*}
from \eqref{ee3}, \eqref{ee1}.  Then in (\ref{ee2}), using (\ref{ee1}) and the fact that the volume of $M$ is equal to 1, we have
$$\int_M e^{-\varphi} \le  \frac{e^{-C_1}}{4} e^{C_1} \int_M e^{-\varphi} + \frac{1}{e} \int_M e^{-\varphi} \le \frac{3}{4} \int_M e^{-\varphi},$$
a contradiction.  Q.E.D.
\end{proof}
\setcounter{remark}{0}
\begin{remark} We note that Lemma \ref{lemmameasure} can also be applied to improve Theorem 1.3 of \cite{TWY}, where it was shown that a conjecture of Donaldson \cite{D2} on \emph{a priori} estimates for the Calabi-Yau equation reduces to estimating a quantity $\int_M e^{-\alpha \varphi}$ for $\alpha>0$, where $\varphi$ is a generalization of the K\"ahler potential function to the non-integrable setting.  Lemma \ref{lemmameasure} can be used to show that it is sufficient to bound, for example, the $L^1$ norm of $\varphi$.
\end{remark}

We now separate the proof of Theorem \ref{tzo} into the two cases (a) and (b), dealing with (b) first.

\bigskip
\noindent
{\it Proof of Theorem \ref{tzo} in the case (b)} \   We have the following lemma.

\begin{lemma} \label{lemmal1}
Suppose that $(M,g)$ is a balanced manifold and $\varphi$ satisfies $\sup_M \varphi=0$ and
$$g'_{i\ov{j}} = g_{i \ov{j}} + \partial_i \partial_{\ov{j}} \varphi >0.$$
Then there exists a uniform constant $C$ depending only on $(M,g)$ such that
\begin{equation}
\int_M | \varphi | d\mu \le C.
\end{equation}
\end{lemma}
\begin{proof}  Observe that this lemma does not require $\varphi$ to solve the complex Monge-Amp\`ere equation (\ref{ma0}).  All we need is that $g$ is balanced and
\begin{equation}
\Delta \varphi > -n.
\end{equation}
Since $g$ is balanced, the canonical Laplacian $\Delta$ coincides (up to a constant multiple) with the Laplace-Beltrami operator for the Riemannian metric associated to $g$ (see for example \cite[Lemma 3.2]{schwarz} or \cite[Proposition 1]{G2}).  Let $G : M \times M \rightarrow \mathbb{R}$ be the associated Green's function.  Then $G$ is uniformly bounded below by $K$, say, and satisfies for all $x$ in $M$,
$$\int_{y \in M} G(x,y) \omega^n(y)=0.$$
Let $x$ in $M$ be a point such that $\varphi(x)=0$.   Then
\begin{eqnarray} \nonumber
0 & = &  \int_M \varphi \, d\mu + \int_{y \in M} G(x,y)(- \Delta \varphi)(y) \omega^n(y)  \\
& < &  \int_M \varphi \, d\mu + nK \int_M \omega^n,  
\end{eqnarray}
and this proves the lemma.  Q.E.D.
\end{proof}

We now finish the proof of Theorem \ref{tzo} in the case (b).  From Lemma \ref{lemmameasure}, we see that 
$$| \varphi| \ge - \inf_M \varphi - C_1 -1$$
on a set of measure at least $e^{-C_1}/4$.  Then 
$$\int_M | \varphi| d\mu \ge \frac{e^{-C_1}}{4} (- \inf_M \varphi - C_1 -1).$$
But $\int_M | \varphi| d\mu$ is bounded from above by Lemma \ref{lemmal1}.  This gives an upper bound for $- \inf_M \varphi$ and completes the proof of Theorem \ref{tzo} in the case (b).

\bigskip
\noindent
{\it Proof of Theorem \ref{tzo} in the case (a)} \  Since we no longer have the balanced condition as in (b), we use a different, and slightly more involved argument, which in fact would also work  in the higher dimensional case.

First recall a basic result of Gauduchon \cite{Ga}.  For any compact Hermitian manifold $(M,\omega)$ of complex dimension $n>1$, there exists a smooth function $u$, unique up to addition of a constant, such that the Hermitian metric $\omega_G=e^u\omega$ satisfies
\begin{equation} \label{gauduchon}
\de\db (\omega_G^{n-1}) =0.
\end{equation}
A metric $\omega_G$ satisfying (\ref{gauduchon}) is called Gauduchon.  We will write $\Delta_G$ for the canonical Laplacian of $\omega_G$.  

We will make use of a Moser iteration result for Gauduchon metrics which can be summarized in the following lemma. Since we are working in the case (a), we restrict now to two complex dimensions.  However, the results given below easily generalize to higher dimensions.

\begin{lemma} \label{lemmapsi} Let $M$ be a two-dimensional compact complex manifold with a Gauduchon metric $\omega_G$.  If $\psi$ is a smooth nonnegative function on $M$ with
$$\Delta_G \psi \ge -C_0$$ 
then there exist constants $C_1$ and $C_2$ depending only on $(M, \omega_G)$ and $C_0$ such that:
\begin{equation} \label{psi1}
\int_M |\partial \psi^{\frac{p+1}{2}}|_{\omega_G}^2 \omega_G^2 \le C_1 p \int_M \psi^p \omega_G^2 \quad \textrm{for all } p \ge 1,
\end{equation}
and
\begin{equation} \label{psi2}
\sup_M \psi \le C_2  \max \left\{  \int_M \psi \, \omega_G^2, 1 \right\}.
\end{equation}
\end{lemma}
\begin{proof}
Compute for $p\ge 1$,
\begin{eqnarray*}
\int_M |\partial \psi^{\frac{p+1}{2}}|_{\omega_G}^2 \omega_G^2 & = &2  \int_M \sqrt{-1} \partial \psi^{\frac{p+1}{2}} \wedge \ov{\partial} \psi^{\frac{p+1}{2}} \wedge \omega_G \\
& = & \frac{(p+1)^2}{2} \int_M \sqrt{-1} \psi^{p-1} \partial \psi \wedge \ov{\partial} \psi \wedge \omega_G \\
& = & \frac{(p+1)^2}{2p} \int_M \sqrt{-1} \partial \psi^p \wedge \ov{\partial} \psi \wedge \omega_G \\
& = & - \frac{(p+1)^2}{2p} \int_M  \psi^p \sqrt{-1} \partial \ov{\partial} \psi \wedge \omega_G + \frac{p+1}{2p} \int_M \sqrt{-1} \, \ov{\partial} \psi^{p+1} \wedge \partial \omega_G \\
& = & \frac{(p+1)^2}{4p} \int_M \psi^p (- \Delta_G \psi) \omega_G^2 \\
& \le &C \frac{(p+1)^2}{4p} \int_M \psi^p \omega_G^2,
\end{eqnarray*}
thus establishing (\ref{psi1}).  Note that, to pass from the fourth to fifth lines, we have made use of the Gauduchon condition (\ref{gauduchon}).

The inequality (\ref{psi2}) then follows by a standard iteration argument, using the Sobolev inequality (\ref{sob})  with $\omega$ replaced by $\omega_G$.
Indeed, write $q=p+1$ and apply (\ref{sob}) to $f=\psi^{q/2}$ to obtain for $q \ge 2$,
$$\left( \int_M \psi^{2q} \omega_G^2 \right)^{1/2} \le C q \max \left\{ \int_M \psi^q \omega_G^2, 1 \right\}.$$
By repeatedly replacing $q$ by $2q$ and iterating we have, after setting $q=2$,
$$ \sup_M \psi \le C \max \left\{ \left( \int_M \psi^{2} \omega_G^2 \right)^{1/2}, 1\right\}.$$
Hence
$$\sup_M \psi \le C \max \left\{  \left( \sup_M \psi \right)^{1/2} \left( \int_M \psi \, \omega_G^2 \right)^{1/2}, 1\right\},$$
and (\ref{psi2}) follows.  Q.E.D.
\end{proof}

We will now apply Lemma \ref{lemmapsi} to the function
\begin{equation}
\psi=\vp-\inf_M\vp \ge 0,
\end{equation}
which, since 
$\Delta\psi>-n$, satisfies
$$\Delta_G \psi = e^{-u} \Delta \psi > -C.$$
Thus, thanks to \eqref{psi2}, to finish the proof of the theorem it suffices to bound the $L^1(\omega_G)$ norm of $\psi$.
From (\ref{psi1}) with $p=1$ we have
\begin{equation}
\int_M |\partial \psi |_{\omega_G}^2 \omega_G^2 \le C  \| \psi \|_{L^1(\omega_G)}.
\end{equation}
Denote by $\underline{\psi}$ the average of $\psi$ with respect to
$\omega_G^2$.  From the Poincar\'e inequality we have
\begin{equation} \label{poincare}
\|\psi-\underline{\psi}\|_{L^2(\omega_G)}\leq C \left( \int_M |\partial \psi |_{\omega_G}^2 \omega_G^2 \right)^{1/2}  \leq C\|\psi\|^{1/2}_{L^1(\omega_G)}.
\end{equation}
We will now make use of Lemma \ref{lemmameasure}.  For a uniform constant $K$, the set 
 $S:=\{  \psi \leq K \} \subset M$ satisfies
 \begin{equation} \label{S}
 |S| \ge \varepsilon,
 \end{equation}
 for a uniform $\ve>0$.
On $S$ we have
$$\underline{\psi}\leq |\psi-\underline{\psi}|+K,$$
and integrating this over $S$ and making use of (\ref{S}), we have
$$\frac{\varepsilon}{\int_M\omega_G^2}\int_M\psi\omega_G^2=\varepsilon \underline{\psi}\leq \int_S \underline{\psi} \omega_G^2 \leq \int_M|\psi-\underline{\psi}| \omega_G^2 +C.$$
Then, using (\ref{poincare}),
\begin{eqnarray*}
\|\psi \|_{L^1(\omega_G)}
 & \leq & C( \|\psi-\underline{\psi}\|_{L^1(\omega_G)}+1) \\
& \leq &   C( \|\psi-\underline{\psi}\|_{L^2(\omega_G)}+1) \\
& \le & C ( \| \psi \|_{L^1(\omega_G)}^{1/2} +1),
\end{eqnarray*}
which shows that $\psi$ is uniformly bounded in $L^1(\omega_G)$.  This  completes the proof of the theorem.  Q.E.D.
\end{proof}

\setcounter{equation}{0}
\section{Proof of the Main Theorem} \label{sectionproof}

From Theorem \ref{tso} and Theorem \ref{tzo}, we have a uniform upper bound of $\tr{g}{g'}$.  It follows from (\ref{ma0}) that $g$ and $g'$ are uniformly equivalent as Hermitian metrics.  To finish the proof of the Main Theorem, it suffices to prove a $C^{\alpha}$ estimate on the metric $g'_{i \ov{j}}$ for some $\alpha>0$.  Indeed once we have this estimate, then differentiating (\ref{ma0}) and applying the standard local elliptic estimates gives us uniform $C^{\infty}$ bounds on $\varphi$.

The $C^{\alpha}$ estimate on $g'_{i\ov{j}}$ is proved in \cite{GL} using a bound on the real Hessian of $\varphi$ (see also \cite{Po,  CKNS, Bl}) and results of Evans \cite{Ev} and Krylov \cite{Kr}.   Instead, we provide a self-contained and more direct proof following a complex version of the Evans-Krylov estimate, due to Trudinger \cite{Tr}.    We follow quite closely Trudinger's argument (see also the expositions of Siu \cite{siubook} and \cite{GT}, and also \cite{W3} for a generalization of this argument to the case of $J$ non integrable).  Since the argument in this case is very similar, we will be brief in places.

We work in a small open subset $U$ of $\mathbb{C}^n$ that contains a ball $B_{2R}$ of radius $2R$, and consider the equation
\begin{equation} \label{ma2}
\log \det (g'_{i \ov{j}}) = \log \det (g_{i \ov{j}} + \varphi_{i \ov{j}}) = F.
\end{equation}
Consider the operator $\Phi$ on positive definite Hermitian matrices given by
$\Phi(A) = \log \det A$
for $A= ( a_{i \ov{j}} )$.  Observe that $\Phi$ is concave.  Our equation (\ref{ma2}) can be expressed as 
\begin{equation} \label{maPhi}
\Phi(g') = F.
\end{equation}
Differentiating  (\ref{maPhi}) with respect to an arbitrary vector $\gamma \in \mathbb{C}^n$ and then $\ov{\gamma}$ we obtain
$$\sum_{i,j} \frac{\partial \Phi  }{\partial a_{i\ov{j}}}g'_{i \ov{j} \gamma} = F_{\gamma}$$
and
$$\sum_{i,j,k,\ell} \frac{\partial^2 \Phi }{\partial a_{i \ov{j}} \partial a_{k \ov{\ell}}} g'_{k \ov{\ell} \ov{\gamma}} g'_{i \ov{j} \gamma} + \sum_{i,j} \frac{\partial \Phi  }{\partial a_{i\ov{j}}} g'_{i \ov{j} \gamma \ov{\gamma}} = F_{\gamma \ov{\ga}}.$$
Note that $$ \frac{\partial \Phi  }{\partial a_{i\ov{j}}} (g') = g'^{i \ov{j}}.$$
Since $\Phi$ is concave, the first term on the left hand side is nonpositive and thus
$$g'^{i \ov{j}} \varphi_{i \ov{j} \gamma \ov{\ga}} \ge F_{\gamma \ov{\ga}} - g'^{i \ov{j}} g_{i \ov{j} \gamma \ov{\ga}}.$$
Writing $w = \varphi_{\gamma \ov{\ga}}$ and 
 $H = F_{\gamma \ov{\ga}} - g'^{i \ov{j}} g_{i \ov{j} \gamma \ov{\ga}}$ 
we have
\begin{equation} \label{w1}
g'^{i \ov{j}} \partial_i \partial_{\ov{j}} w \ge H.
\end{equation}
Note that for each fixed $\gamma$, $H$ is uniformly bounded.  

Before we make use of (\ref{w1}), we observe that another consequence of the concavity of $\Phi$ is that for all $x,y$ in $U$, 
$$\Phi(g'(y)) + \sum_{i,j} \frac{\partial \Phi}{\partial a_{i \ov{j}}}(g'(y)) \left( g'_{i \ov{j}}(x) - g'_{i \ov{j}}(y) \right) \ge \Phi(g'(x)).$$
It follows that
\begin{equation} \label{concavity2}
\sum_{i,j} \frac{\partial \Phi}{\partial a_{i \ov{j}}}(g'(y))  \left( g'_{i \ov{j}}(y) - g'_{i \ov{j}}(x) \right) \le F(y)-F(x) \le CR,
\end{equation}
where for the last line we have used the mean value theorem.
Since we have \emph{a priori} estimates on $g'$ and hence on $\partial \Phi / \partial a_{i \ov{j}}$, we can find  unit vectors $\gamma_1, \ldots, \gamma_N$ in $\mathbb{C}^n$ and real-valued functions $\beta_1, \ldots, \beta_N$ satisfying
$$0<\frac{1}{C} \le \beta_\nu \le C, \quad \textrm{for } \nu =1, \ldots, N,$$
such that 
\begin{equation}
\frac{\partial \Phi}{\partial a_{i \ov{j}}}(g'(y)) = \sum_{{\nu}=1}^N \beta_{\nu}(y) (\gamma_{\nu})^i \ov{(\gamma_{\nu})^j}.
\end{equation}
Moreover, we may assume that $\gamma_1, \ldots, \gamma_N$ contains an orthonormal basis for $\mathbb{C}^n$.
From (\ref{concavity2}) and the mean value theorem again we obtain
\begin{equation} \label{concavity3}
\sum_{{\nu}=1}^N \beta_\nu (w_{\nu}(y) - w_\nu(x)) \le CR, \quad \textrm{for} \quad w_{\nu} = \varphi_{\gamma_{\nu} \ov{\gamma_{\nu}}}.
\end{equation}

We will obtain the desired estimate from (\ref{w1}), (\ref{concavity3}) and 
 the following Harnack inequality (for its proof see Theorem 9.22 of \cite{GT}).

\begin{lemma} \label{harnack}
Let $g'$ be a Hermitian metric on $U \subset \mathbb{C}^n$ which is uniformly equivalent to the Euclidean metric.  Suppose that $v \ge 0$ satisfies
$$g'^{i \ov{j}} \partial_i \partial_{\ov{j}} v \le \theta,$$
on $B_{2R} \subset U$. Then  there exist uniform constants $p>0$ and $C>0$ such that
$$\left( \frac{1}{R^{2n}} \int_{B_R} v^p \right)^{1/p} \le C \left( \inf_{B_R} v + R \| \theta \|_{L^{2n}(B_{2R})} \right).$$
\end{lemma}

For $j=1,2$, write $$M_{j\nu} = \sup_{B_{jR}} w_{\nu}, \quad m_{j\nu} = \inf_{B_{jR}} w_{\nu}, \quad \textrm{and} \quad \Omega(jR) = \sum_{\nu=1}^N \textrm{osc}_{B_{jR}} w_{\nu} =\sum_{\nu=1}^N  (M_{j\nu}- m_{j\nu}).$$  Since each $w=w_{\nu}$ satisfies (\ref{w1}), we can apply 
Lemma \ref{harnack} to $M_{2\nu}-w_{\nu}$ to obtain
\begin{equation} \label{Mw}
\left(\frac{1}{R^{2n}}\int_{B_R}(M_{2\nu}-w_{\nu})^p \right)^{1/p} \le C\left(M_{2\nu}-M_{1\nu} + R^2\right).
\end{equation}
Thus for a fixed $\ell \in \{1, \ldots, N\}$, we have
\begin{eqnarray} \nonumber
\left( \frac{1}{R^{2n}}\int_{B_R}(\sum_{\nu \neq \ell} (M_{2\nu}-w_{\nu}))^p \right)^{1/p} & \le & N^{1/p} \sum_{\nu \neq \ell} \left(\frac{1}{R^{2n}}\int_{B_R}(M_{2\nu}-w_{\nu})^p \right)^{1/p} \\ \nonumber
& \le & C \left( \sum_{\nu \neq \ell} (M_{2\nu} - M_{1\nu}) + R^2 \right) \\ \label{sum}
& \le & C \left( \Omega(2R)- \Omega(R) + R^2 \right),
\end{eqnarray}
since
$$(M_{2\nu} - m_{2\nu}) - (M_{1\nu} - m_{1\nu}) = (M_{2\nu}-M_{1\nu}) + (m_{1\nu}-m_{2\nu}) \ge M_{2\nu} - M_{1\nu}.$$
But from (\ref{concavity3}), we have
$$\beta_\ell (w_\ell (y)-w_\ell(x)) \le CR + \sum_{\nu \neq 
\ell} \beta_{\nu} (w_{\nu}(x)-w_{\nu}(y)).$$
If we choose $x$ so that $w_\ell(x)$ approaches $m_{2\ell}$ we obtain
$$w_\ell(y) - m_{2\ell} \le C\left(R + \sum_{\nu \neq \ell} (M_{2\nu} - w_{\nu}(y))\right).$$
Integrating in $y$ over $B_R$ and applying (\ref{sum}) gives
\begin{equation} \label{mw}
\left( \frac{1}{R^{2n}} \int_{B_R} (w_\ell-m_{2\ell})^p \right)^{1/p} \le C\left( \Omega(2R)- \Omega(R) +R \right).
\end{equation}
On the other hand, from (\ref{Mw}) we have
\begin{equation} \label{Mw2}
\left(\frac{1}{R^{2n}}\int_{B_R}(M_{2\ell}-w_{\ell})^p \right)^{1/p} \le C\left(\Omega(2R)- \Omega(R) +R\right).
\end{equation}
Adding (\ref{mw}) and (\ref{Mw2}) and summing over $\ell$ we have
$$\Omega(2R) \le C\left( \Omega(2R)- \Omega(R) + R \right),$$
and hence
$$\Omega(R) \le \delta \Omega(2R) + CR$$
for a uniform $0< \delta<1$.  It then follows by a standard argument (see \cite{GT}, Chapter 8) that there exist uniform constants $C$ and $\kappa>0$ such that 
$$\Omega(R) \le C R^{\kappa}.$$
The desired H\"older estimate on $\varphi_{i \ov{j}}$ and hence $g'$ follows.  This completes the proof of the main theorem.

\setcounter{remark}{0}

\begin{remark} In \cite{ChP}, Cherrier used the maximum principle to prove an analogue of Yau's \cite{Y} third order estimate of $\varphi$  for  Hermitian $\omega$.  For some related results, we refer the reader to:  Calabi's paper  \cite{Ca2} which inspired Yau's computation;  \cite{PSS} for a concise proof of the parabolic version of this estimate; \cite{adiab} for a precise version of this estimate without assuming a $C^2$ estimate;
and \cite{TWY} for a generalization of the third order estimate in the setting of Donaldson's program \cite{D2}.

\end{remark}

\setcounter{remark}{0}
\setcounter{equation}{0}
\section{Solving the Monge-Amp\`ere equation} \label{sectionopenness}

In this section we give a proof of Corollary 1.  We use the continuity method and consider the family of equations
\begin{equation}\label{cty1}
(\omega + \sqrt{-1} \partial \ov{\partial} \vp_t)^n = e^{tF+b_t} \omega^n, \quad \textrm{for } t \in [0,1],
\end{equation}
with 
\begin{equation} \label{pos}
\omega+ \sqrt{-1} \partial \ov{\partial} \varphi_t>0,
\end{equation} 
where $b_t$ is a constant for each $t$.  Fix $\alpha$ with $0<\alpha<1$.
We claim that  (\ref{cty1}),  (\ref{pos}) can be solved for  $\varphi_t \in C^{2, \alpha}(M)$ and $b_t \in \mathbb{R}$, for all $t \in [0,1]$.  
Consider the set
\begin{eqnarray*}
\mathcal{T} & =&  \{ t' \in [0,1] \ | \ \textrm{there exists } \varphi_t \in C^{2, \alpha}(M) \ \textrm{and } b_t   \\
&& \mbox{}  \ \quad \qquad \qquad \textrm{solving (\ref{cty1}) and (\ref{pos}) for } t\in [0,t'] \}.
\end{eqnarray*}

Clearly, $0 \in \mathcal{T}$.  If we can show that $\mathcal{T}$ is both open and closed in $[0,1]$ then, using the higher order estimates given above, this proves the theorem.  
We first show that $\mathcal{T}$ is open.  Assume $\hat{t}$ is in $\mathcal{T}$.  We will show that there exists $\varepsilon>0$ such that $t \in \mathcal{T}$ for $t \in [\hat{t},\hat{t}+\varepsilon)$.

Write $\hat{\omega}$ for the Hermitian metric  $$\hat{\omega}= \omega + \sqrt{-1} \partial \ov{\partial} \varphi_{\hat{t}}>0.$$  We wish to show that there exists  $\psi_t \in C^{2,\alpha}$ for $t \in [\hat{t}, \hat{t} + \varepsilon)$  with $\psi_{\hat{t}}=0$, solving
$$(\hat{\omega} + \sqrt{-1} \partial \ov{\partial} \psi_t)^n = e^{(t-\hat{t})F+ b_t-b_{\hat{t}}} \hat{\omega}^n,$$
for $b_t$ a  function of $t$.

Applying Gauduchon's theorem to $\hat{\omega}$, we see that there is a function $\hat{u}$ such that 
$\hat{\omega}_G = e^{\hat{u}} \hat{\omega}$ is Gauduchon.  We may assume by adding a constant to $\hat{u}$ that 
\begin{equation} \label{ucondition}
\int_M e^{(n-1)\hat{u}} \hat{\omega}^n =1.
\end{equation} 
We will show that we can find $\psi_t \in C^{2,\alpha}(M)$ for $t \in [\hat{t}, \hat{t}+\varepsilon)$ solving
\begin{equation} \label{equation1}
(\hat{\omega} +\sqrt{-1} \partial \ov{\partial} \psi_t)^n = \left( \int_M e^{(n-1)\hat{u}}(\hat{\omega}+ \sqrt{-1} \partial \ov{\partial} \psi_t)^n \right) e^{(t-\hat{t})F +c_t} \hat{\omega}^n, \quad \textrm{for } t \in [\hat{t},\hat{t}+\varepsilon), 
\end{equation}
where $c_t$ is chosen so that  
$$\int_M e^{(t-\hat{t})F+c_t} e^{(n-1)\hat{u}} \hat{\omega}^n=1.$$
Notice that  $c_{\hat{t}}=0$.

Define two Banach manifolds $B_1$ and $B_2$ by
$$B_1 = \{ \psi \in C^{2,\alpha}(M) \ | \ \int_M \psi e^{(n-1)\hat{u}} \hat{\omega}^n =0 \}, \quad B_2 = \{ h \in C^{\alpha}(M) \ | \ \int_M e^{h} e^{(n-1)\hat{u}} \hat{\omega}^n = 1 \}.$$
Define a linear operator $\Psi: B_1 \rightarrow B_2$ by
$$\Psi(\psi) = \log \frac{(\hat{\omega} + \sqrt{-1} \partial \ov{\partial} \psi)^n}{\hat{\omega}^n} - \log \left( \int_M e^{(n-1)\hat{u}} (\hat{\omega} + \sqrt{-1} \partial \ov{\partial} \psi)^n \right).$$
We wish to find $\psi_t$ solving $\Psi(\psi_t) = (t-\hat{t})F + c_t$ for $t \in [\hat{t}, \hat{t} + \varepsilon)$.

Observe that $\Psi(0)=0$.  By the inverse function theorem, to get openness of the equation (\ref{equation1}) at $\psi=0$ we just need to show that $$(D \Psi)_0 : T_0 B_1=B_1  \longrightarrow T_0 B_2 = \{ \rho \in C^{\alpha}(M) \ | \ \int_M \rho e^{(n-1)\hat{u}} \hat{\omega}^n=0 \}$$ is invertible.  Compute
$$(D \Psi)_0 (\eta ) = \Delta_{\hat{\omega}} \eta - n\int_M e^{(n-1)\hat{u}} \hat{\omega}^{n-1} \wedge \sqrt{-1} \partial \ov{\partial} \eta = \Delta_{\hat{\omega}} \eta,$$
where we are using the fact that $$\partial \ov{\partial} ( e^{(n-1) \hat{u}} \hat{\omega}^{n-1} ) = \partial \ov{\partial} (\hat{\omega}_{\textrm{G}}^{n-1})=0.$$
Recall that we can solve the equation $\Delta_{\hat{\omega}_{\textrm{G}}} f = v$ as long as $\int_M v \hat{\omega}_{\textrm{G}}^n=0$ (see e.g. \cite{Bu}). Given $\rho \in T_0 B_2$ we see that 
$$\int_M \rho e^{-\hat{u}} \hat{\omega}_{\textrm{G}}^n =0.$$
Then if $\eta$ solves $$e^{-\hat{u}} \Delta_{\hat{\omega}} \eta = \Delta_{\hat{\omega}_{\textrm{G}}} \eta = \rho e^{-\hat{u}}$$ we have
$$\Delta_{\hat{\omega}} \eta = \rho,$$
as required.  This shows that we can find $\psi_t$ solving (\ref{equation1}) for $t \in [\hat{t}, \hat{t}+\varepsilon)$ as required.  Thus $\mathcal{T}$ is open.

To prove that $\mathcal{T}$ is closed we need estimates on both $\varphi_t$ and $b_t$.  Because of our Main Theorem,  it suffices to show that $b_t$ is uniformly bounded.  But from (\ref{cty1}) we see that if $\varphi_t$ achieves a maximum at a point $x$ in $M$ we have 
$tF(x) +b_t \le 0$.  Combining this with a similar argument at  minimum  of $\varphi_t$, we obtain
$$|b_t| \le \sup_M |F|.$$
This completes the proof of Corollary 1.

\begin{remark} \label{unique}
We end this section with a remark about uniqueness, that dates back to \cite{C}.  Solutions $\varphi$ to (\ref{ma0}), (\ref{sup}) are unique for a given fixed $F$.  To see this, let $\varphi_1$ and $\varphi_2$ satisfy
\begin{equation}
(\omega+\sqrt{-1} \partial \ov{\partial} \varphi_1)^n = (\omega+ \sqrt{-1} \partial \ov{\partial} \varphi_2)^n,
\end{equation}
with $\sup_M \varphi_1=0=\sup_M \varphi_2$.
Then, writing $\omega_j= \omega+ \sqrt{-1}\partial \ov{\partial} \varphi_j$ for $j=1,2$ and $\theta = \varphi_1-\varphi_2$, we see that
$$\sqrt{-1} \partial \ov{\partial} \theta \wedge \sum_{i=0}^{n-1} \omega_1^i \wedge \omega_2^{n-1-i}=0,$$
and it follows by the strong maximum principle that $\theta=0$.   Moreover, the constant $b$ in (\ref{ma3}) is unique by a simple maximum principle argument (see \cite{ChP}).
\end{remark}

\setcounter{equation}{0}
\section{Prescribing the first Chern form} \label{sectionchern}

In this section, we prove Corollary 2.  First we describe some terminology. If $g$ is a Hermitian metric and $\omega$ its associated $(1,1)$-form, we will denote by $\Ric(\omega)$ the first Chern form of the canonical connection of $g$, defined as follows: if $R^i_{jk\ov{\ell}}$ is the curvature of the canonical connection of $g$, we let $R_{k\ov{\ell}}=\sum_i R^i_{ik\ov{\ell}}$ and
$$\Ric(\omega)=\frac{\sqrt{-1}}{2\pi}\sum_{k,\ell} R_{k\ov{\ell}}dz^k\wedge d\ov{z^\ell}.$$
$\Ric(\omega)$ is a closed real $(1,1)$-form, cohomologous to the first Chern class 
$c_1(M)$.

We now give the proof of Corollary 2.  First notice the standard transgression formula (see e.g. \cite[(3.16)]{TWY})
\begin{equation}\label{bott}
\Ric(\omega+\sqrt{-1}\de\db\varphi)-\Ric(\omega)=
-\frac{\sqrt{-1}}{2\pi}\de\db\log\frac{(\omega+\sqrt{-1}\de\db\varphi)^2}{\omega^2},
\end{equation}
and so integration by parts shows that \eqref{eqn} implies \eqref{constr}.

Conversely, assume that \eqref{constr} holds. Let $\Delta_G$ be the canonical Laplacian of $\omega_G$. Because of \eqref{constr}, there exists a smooth function $f$ such that 
$$\frac{\Delta_G f}{2\pi}=\frac{2(\Ric(\omega)-\psi)\wedge\omega_G}{\omega_G^2}.$$
We claim that in fact we have 
\begin{equation}\label{chern}
\Ric(\omega)=\psi+\frac{\sqrt{-1}}{2\pi}\de\db f.
\end{equation}
To see this, call $a=\Ric(\omega)-\psi-\frac{\sqrt{-1}}{2\pi}\de\db f$, and notice that $a$ is an exact real $(1,1)$ form that satisfies
$$\omega_G\wedge a=0.$$
Thus $a$ is $\omega_G$-anti-selfdual, and since it is exact it has to be zero because $$0=\int_M a^2=-\|a\|^2_{L^2(\omega_G)}.$$
By Corollary 1 there exist a constant $b$ and a Hermitian metric $\omega+\sqrt{-1}\de\db\varphi$ such that
\begin{equation}\label{ma4}
(\omega+\sqrt{-1}\de\db\varphi)^2=e^{f+b}\omega^2.
\end{equation}
Combining \eqref{ma4} with \eqref{bott} and \eqref{chern} gives \eqref{eqn}. This completes the proof of Corollary 2.

\setcounter{remark}{0}
\begin{remark}
It follows from the proof of Corollary 2 that \eqref{constr} is equivalent to the statement that $\Ric(\omega)-\psi$ is $\de\db$-exact, and thus if \eqref{constr} holds for one Gauduchon metric then it holds for all.
Clearly, in the K\"ahler case \eqref{constr} is always satisfied.
\end{remark}

\bigskip
\noindent
{\bf Acknowledgements} \ 
The authors would like to thank S.-T. Yau for many useful discussions on the complex Monge-Amp\`ere equation and also on the balanced condition. The authors express their gratitude to D.H. Phong for his support, encouragement and helpful suggestions.  In addition the authors thank V. Apostolov, L. Ni, J. Song, G. Sz\'ekelyhidi and L.-S. Tseng for some helpful discussions.  We thank X. Zhang for sending us his  related preprint \cite{Zh},  from which we also learned about the reference \cite{ChP}.

\bigskip
\noindent
Mathematics Department, Columbia University, 2990 Broadway, New York, NY 10027

\bigskip
\noindent
Mathematics Department, University of California, San Diego, 9500 Gilman Drive \#0112, La Jolla CA 92093

\end{document}